\documentclass{article}

\usepackage{times}

\usepackage{algorithm}
\usepackage{algorithmic} 
\usepackage{amsmath}
\usepackage{amssymb}
\usepackage{graphicx}
\usepackage{subfigure}
\newtheorem{prop}{Proposition.}

\title{Accelerated Kaczmarz Algorithms using History Information}
\author{ \\ Tengfei Ma
  \\ IBM Research-Tokyo
  \\feitengma0123@gmail.com
}

\begin{document}

\maketitle

\begin{abstract}
 The Kaczmarz algorithm is a well known iterative method for solving overdetermined linear systems. Its randomized version yields provably exponential convergence in expectation. In this paper, we propose two new methods to speed up the randomized Kaczmarz algorithm by utilizing the past estimates in the iterations. The first one utilize the past estimates to get a preconditioner. The second one combines the stochastic average gradient (SAG) method with the randomized Kaczmarz    algorithm. It takes advantage of past gradients to improve the convergence speed. Numerical experiments
indicate that the new algorithms can dramatically outperform the standard randomized Kaczmarz algorithm. 
\end{abstract}

\section{Introduction}
 The Kaczmarz algorithm (\cite{kaczmarz1937angenaherte}) is a simple but powerful iterative method for solving the overdetermined system with equations $Ax=b$. Due to its simplicity and speed, it has a wide range of applications from computer tomography to image reconstruction (\cite{sezan1987applications}). It is a form of alternating projection method, which in each iteration projects the current solution to a subspace.

Given a matrix $A \in \mathbb{R}^{m\times n}$ with $m \ge n$ and $b \in \mathbb{R}^m$, we denote the rows of $A$ by $a_1^T,a_2^T,..a_m^T$ and $b = (b_1,b_2,...,b_m)^T$. The Kaczmarz method project the current estimation orthogonally onto the solution hyperplane of $a_j^T x =b_j$, where the row $j$ is selected in a cyclic manner.

Recently, \cite{strohmer2009randomized} proposed to select the row with biased sampling and proved that the randomized Kaczmarz method (RK) converges with expected exponential rate. Let $||A||_F^2$ denote the Frobenius norm of $A$ and $||\cdot||$ denote the standard norm. In each iteration a row $i$ is randomly selected with probability proportional to $||a_i||^2$, and finally we could get the following exponential bound for the convergence in expectation:
\begin{equation}
\label{converRate}
\mathbb{E}||x_k-x||^2 \le (1-\frac{1}{\kappa(A)^2})^k ||x_0-x||^2
\end{equation}
where $\kappa (A) = ||A||_F||A^{-1}||$, $A^{-1}$ is the left inverse of $A$ which is always assumed to exist, and $x_0$ is an arbitrary initial value. Using this algorithm, the cost per iteration is $O(n)$ and the expected iteration for convergence is $O(\log(1/\epsilon))$ where $\epsilon$ is the accuracy parameter. So it is computationally feasible for very large systems. It is also shown that the randomized Kaczmarz method often outperforms the celebrated conjugate gradient method.

Besides consistent linear system, the RK algorithm has also been analyzed for inconsistent linear system $Ax = b + w$ where $w$ is an arbitrary noise vector (\cite{needell2010randomized}). And some extended RK methods are proposedd for systems of linear inequalities (\cite{leventhal2010randomized}), least square problems (\cite{zouzias2013randomized}), and online compressed sensing (\cite{lorenz2014sparse}).
The RK algorithm has a theoretical linear convergence rate. However, the convergence rate largely depends on the condition number $\kappa$ of matrix $A$, and the convergence will be extremely slow for ill-conditioned problems. Therefore, some accelerated RK methods are proposed. For example, \cite{liu2015accelerated} applied the Nesterov acceleration scheme to the standard RK algorithm, and obtained the accelerated randomized Kaczmarz algorithm (ARK). 

In this paper, we develop new acceleration schemes for the Kaczmarz algorithm by utilizing history information. The basic idea is to change the direction of projection in each iteration to make it converge faster. In the RK algorithm, each projection is always along a row vector $a_j$ which is orthogonal to the hyperplane $a_j^T x=b_j$. Our first acceleration scheme finds a new preconditioner $C$ which changes the projection direction into $C a_j$. The preconditioner is approximated based on the estimate of $x$ in past iterations. Our second acceleration scheme considers the relationship between stochastic gradient descent (SGD) and the randomized Kaczmarz algorithm and combines them, In each iteration, we first use the past gradients to get a variant of SGD, the stochastic average gradient (SAG). We do a gradient descent along the SAG and then project the point back into a hyperplane.

The paper is organized as follows. The next section covers related work about accelerated Kaczmarz algorithms. In section 3, we introduce the preconditioning technique and induce our new preconditioner. Section 4 present the second acceleration scheme which integrates SAG into the RK algorithm. Numeric experiments are shown in Section 5 and we conclude the paper in Section 6.

\begin{algorithm}
\caption{Randomized Kaczmarz Algorithm}
\label{RKalg}
\begin{algorithmic}[1]
\STATE Initialize $k \leftarrow 0$
\FOR {$k=0,1,...$}
\STATE Select row $j$ from $\{1,2,...m\}$ with probability $\frac{||a_j||^2}{||A||_F^2}$
\STATE Project $x_{k+1} = x_k + \frac{(b_j - a_j^T x_k)}{||a_j||^2} a_j$
\STATE Update $k \leftarrow k+1$
\ENDFOR
\end{algorithmic}
\end{algorithm}

\section{Related Work}

\subsection{Improvement to the Kaczmarz algorithm}
Since the RK was analyzed by \cite{strohmer2009randomized}, there has been several directions to extend it. Two-subspace projection method extends RK by iterately projecting the estimate onto the solution space given by two randomly selected rows. It only improves when the system has correlated rows. \cite{eldar2011acceleration} has a different strategy to select the rows. It projects the row vectors onto a low dimensional space, and then selects the row which leads to the largest improvement. Beyond the scope of randomized Kaczmarz algorithm, there is some work focusing on accelerating the classical Kaczmarz method. For example, Brezinski and Redivo–Zaglia \cite{brezinski2013convergence} use sequence transformation to change the projection procedure. But it is too complex to get a transformed sequence and it needs to store too many vectors, so it is difficult to be applied in practice when the dimension of $A$ is large. 

\subsection{The Randomized Kaczmarz Algorithm and Stochastic Gradient Descent}
The Kaczmarz algorithm is fundamentally a special case of alternating projection (\cite{strohmer2009randomized}). In some area it is called POCS (projection to convex sets).  But it also has a strong relationship with stochastic gradient descent (SGD). Very recently, It has been demonstrated that the RK algorithm is equivalent to a form of stochastic gradient descent with weighted sampling (\cite{needell2014stochastic}). However, one advantage of the RK is that it does not need to set up the step size for each iteration, although it generally does not get the optimal step size.

Considering the connection between the RK and the SGD, Liu and Stephen \cite{liu2015accelerated} apply the well known Nesterov's acceleration procedure to the RK algorithm. They demonstrate the convergence of their accelerated randomized kaczmarz algorithm (ARK) and obtain significant improvement for ill-conditioning problems in numeric experiments (\cite{liu2015accelerated}). The ARK introduces two additional sequences $\{y_k\}$ and $\{v_k\}$ as follows
\begin{eqnarray*}
y_k &=& \alpha_k v_k + (1-\alpha_k) x_k\\
x_{k+1} &=& y_k - a_i(a_i^T y_k-b_i)/||a_i||^2\\
v_{k+1} &=& \beta_k v_k + (1-\beta_k) y_k - \gamma a_i (a_i^T y_k-b_i)/||a_i||^2
\end{eqnarray*}
where the scalars $\alpha_k$,$\beta_k$ and $\gamma_k$ are calculated offline based on the hyperparameter $\lambda \in [0,\lambda_{min}]$, ($\lambda_{min}$ is the minimum eigenvalue of $A^T A$). They prove that when $\lambda >0$, the ARK gets a linear convergence rate. The ARK is then extended to solving the sparse data. 

The ARK performs very well for ill-conditioned problems, especially when the $\lambda_{min}(A^T A)$ is known. However, to get the accurate  $\lambda_{min}(A^T A)$ is difficult. And in many cases the inaccurate $\lambda_{min}(A^T A)$ will lead to much worse performance.

\subsection{Importance of History Information}
History information has been used in many acceleration schemes for alternating projection (\cite{gearhart1989acceleration}) and SGD (\cite{roux2012stochastic}, \cite{johnson2013accelerating}, \cite{nitanda2014stochastic}). The main idea to use history information is to find a point closest to the final solution in each iteration (\cite{gearhart1989acceleration}) or reduce the variance between stochastic gradients and the full gradients (\cite{johnson2013accelerating}). The ARK also utilizes the history information by employing a Nesterov's acceleration procedure. 
In this paper, we developed two approaches to utilizing history information. We use the past estimates of $x$ to approximate a preconditioner in our first algorithm, while in the second algorithm we use the past stochastic gradients to approximate full gradients as in (\cite{roux2012stochastic}).

\section{Approximated Preconditioned Kaczmarz(APK) Algorithm}
\label{sec:APK}
The motivation of our first acceleration scheme lies on two aspects. Firstly, we consider using preconditioning to reduce the condition number of $A$ in the system. Secondly, we want to use the history information to generate a proper preconditioner.

As we explained before, the convergence rate of the randomized Kaczmarz algorithm largely depends on the condition number of $A$. When this number is large, the convergence speed will be too slow. One solution to this problem is to use a preconditioning matrix. A preconditioning matrix (or preconditioner) $B$ of a matrix $A$ is a matrix such that $BA$ or $AB$ has a smaller condition number than $A$. So the original problem could be changed into either a left preconditioned system $$BAx = Bb$$ or a right preconditioned system.
\begin{equation}
\label{eq:RPrecondition}
ABB^{-1}x = b
\end{equation}.

Here we consider the right preconditioned system \ref{eq:RPrecondition}. Assume that we already know the preconditioner $B \in \mathbb{R}^{n*n}$, we explain each iteration of the new Kaczmarz algorithm as follows.

First we solve the new linear system $ABy = b$, where $y = B^{-1}x$. So at each iteration we project the current estimation $y_k$ on the hyperplane defined by the row $i$: 
\begin{equation}
y_{k+1} = y_k + \frac{b_i-a_i^T B y_k}{a_i^T BB^T a_i} B^T a_i
\end{equation}.
Replace $y_k$ with $ B^{-1}x_k$, then we get the update of $x$:
\begin{eqnarray}
B^{-1} x_{k+1} &=& B^{-1} x_k + \frac{b_i-a_i^T B B^{-1} x_k}{a_i^T BB^T a_i} B^T a_i\\
x_{k+1} &=& x_k + \frac{b_i-a_i^T x_k}{a_i^T BB^T a_i} BB^T a_i
\end{eqnarray}.

A good preconditioner may accelerate the convergence a lot. Indeed, the choice of preconditioner is often more important than the choice of iterative method, according to Yousef Saad (\cite{saad2003iterative}). However, how to select a preconditioning matrix remains a difficult problem. In many cases determining a good preconditioning matrix itself has the same computational complexity with the original problem. 

In order to use a preconditioner in the Kaczmarz algorithm, it is better to keep the preconditioner to be a diagonal matrix. Each iteration of the Kaczmarz algorithm costs only $O(n)$. If the preconditioner $B$ is not diagonal or sparse enough, the computation of $a_i^T B$ costs $O(n^2)$, which will be unacceptable for just one iteration. Another choice is to directly get $ BB^T a_i$ without an previous estimation of $B$, such as the online LBFS (oLBFS) method (\cite{schraudolph2007stochastic}), which is a stochastic quasi-newton method. But it still costs more than using a diagonal matrix.

We aim to get a diagonal matrix $C = BB^T$, thus the ''projection'' could be represented by a modified form which only contains $C$:
\begin{equation}
x_{k+1} = x_k + \frac{b_i-a_i^T x_k}{a_i^T C a_i} C a_i
\label{update}
\end{equation}
Left-multiply the two sides with a vector $a_i^T$, we find that $a_i^T x_{k+1} = b_i$. That means $x_{k+1}$ is still on the hyperplane given by the original row $i$.

\subsection{Optimizing the diagonal preconditioner using history information}
Consider that the Kaczmarz algorithm is essentially an alternating projection method. After each iteration, the new estimation lay on a hyperplane. 


Assume that we have a list of past estimates in the classical Kaczmarz algorithm $x_1,...x_m, x_{m+1}, x_{2m}$, where rows are selected with a cyclic manner according to an order $R(1,...m)$. So each $x_k$ is projected onto the hyperplane given by the row $i_{k+1} = R(k+1)$ and leads to the next estimation $x_{k+1}$. Since we do not change the selection order, $x_k$ and $x_{k+m}$ are on the same hyperplane. 

\begin{figure}[htb]
\centering
\includegraphics[width=\linewidth]{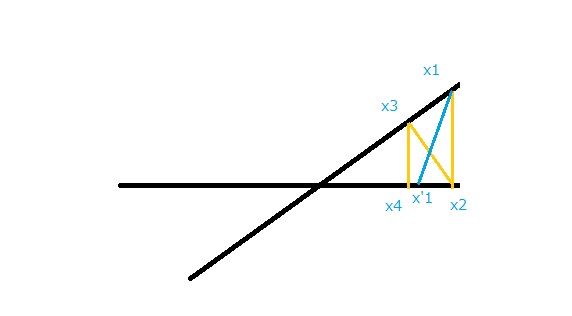}
\caption{A simple example of the preconditioning idea. In this simple $\mathbb{R}^{2*2}$ case, we want to use a new projection vector. It projects $x_1$ to $x'_1$ which is close to $x_4$}
\label{fig:example}
\end{figure}

The idea of our method is, why do not we directly use a pseudo projection to project $x_k$ to $x_{k+m+1}$ instead of the original $x_{k+1}$? That means, we let the pseudo projection seemingly jump across a cycle of projections. So we use the preconditioned row vector $a_{i_k} C$ as the direction of the pseudo projection. we show a very simple case in Figure~\ref{fig:example} as an example.
\begin{equation}
\label{eq:obj0}
C = \arg \min_C F(C) =  \arg \min \sum_{k=2}^{m} ||x_{k+m}-x'_{k}||^2
\end{equation}
where $x'_{k}$ is a projection of $x_{k-1}$ along the direction of $a_{i_k}C$ instead of $a_{i_k}$. So the objective function becomes
\begin{equation}
\label{eq:obj}
F_1(C) = \sum_{k=2}^{m} ||x_{k+m}-x_{k-1} -  \frac{b_{i_k}-a_{i_k}^T  x_{k-1}}{a_{i_k}^T C a_{i_k}} C a_{i_k} ||^2
\end{equation}
To simplify the optimization problem, we have the following strategy to approximate the objective function. We assume that $C$ is not too distant from $I$, so that we could have $a_i^T C a_i \simeq a_{i_k}^T a_{i_k} $. In this case, the objection function has been changed into a combination of two parts, an approximation of the Equation\ref{eq:obj}, and a regularization term to keep $a_{i_k}^T C a_{i_k}$ close to $a_{i_k}^T a_{i_k}$. To keep similarity, we used the Frobenius norm of (C-I) as regularization.
\begin{equation}
\label{eq:FC}
 \sum_{k=2}^{m} ||x_{k+m}-x_{k-1} -  \frac{b_{i_k}-a_{i_k}^T  x_{k-1}}{a_{i_k}^T a_{i_k}} C a_{i_k} ||^2 + \alpha||C-I||_F^2
\end{equation}
where $\alpha$ is a coefficient parameter for the regularization. 

As $C$ is diagonal, we only need to calculate the diagonal vector of $C$. Denote $s = diag(C)$ as the diagonal vector of $C$, and $A_{i_k}$ as a diagonal matrix whose diagonal is $a_{i_k}$. Then $C a_{i_k}$ could be transformed into $A_{i_k} s$, and $F_2(C)$ can be written as functions of $s$: $F_2(C) = F(s)$.

As the objective function turned to be convex, it is easy to get the solution by making the derivative $F'(s) = 0$. From the objective function 
\begin{equation}
F(s) = \sum_{k=2}^{m} ||x_{k+m}-x_{k-1} -  \frac{b_{i_k}-a_{i_k}^T  x_{k-1}}{a_{i_k}^T  a_{i_k}}  A_{i_k} s ||^2 + \alpha || (s - e) ||^2
\end{equation}
where $e$ is a all-ones vector, we differentiate $F(s)$ with respect to $s$, and get the derivative as
\begin{equation}
F'(s) =  F'_1(s) + 2 \alpha (s - e)
\end{equation}
where 
\begin{equation}
F'_1(s) = -2 \sum_{k=2}^{m} (\frac{b_{i_k}-a_{i_k}^T x_{k-1}}{a_{i_k}^T  a_{i_k}}) A_{i_k} (x_{k+m}-x_{i_k} -  \frac{b_{i_k}-a_{i_k}^T  x_{k-1}}{a_{i_k}^T  a_{i_k}}  A_{i_k} s)
\end{equation}
Let $F'(s) = 0$, $A'_{i_k} = \frac{b_{i_k}-a_{i_k}^T  x_{k-1}}{a_{i_k}^T  a_{i_k}}  A_{i_k}$, and $\delta_{i_k} = x_{k+m}-x_{k-1}$, we get the optimal $s$:
\begin{eqnarray}
s &=& \arg \min_s F_1(s) \\\nonumber
 &=& \left(\sum_{k=2}^{m} (\alpha + A_{i_k}^{'2})\right)^{-1} \left(\sum_{k=2}^{m} A'_{i_k} \delta_{i_k} + \alpha e \right)
\end{eqnarray}
Note that $A'_{i_k}$ is a diagonal matrix, so the computation costs only $O(mn)$. And when $m$ is extremely large, we can use only a subset of samples from {1,...m} instead of a full computation. In practice we update this diagonal matrix after a long interval, so the costs does not impact a lot.

\subsection{Convergence analysis of the APK Algorithm}
\begin{prop}
The APK algorithm converges at a linear rate.
\end{prop}
The conclusion is very intuitive. It can be easily proved from the convergence analysis of the randomized Kaczmarz algorithm \ref{converRate}. We apply the conclusion (\ref{converRate}) to the right reconditioned form (\ref{eq:RPrecondition}) which we use in our APK algorithm: 
\begin{equation}
\mathbb{E}||B^{-1} x_k- B^{-1} x||^2 \le (1-\frac{1}{\kappa(AB)^2})^{k-t0} ||x_{t0}-x||^2
\end{equation}
where $x_{t0}$ is the start point after we update the preconditioner $C$, and $x_k$ is an estimate before we update $C$ next time.
Denote the minimum and maximum eigenvalues of $C = B*B$ as $\lambda_{min}(C)$ and $\lambda_{max}(C)$ separately. Then $||B^{-1} x_k- B^{-1} x||^2 \ge \lambda_{max}(C)^{-1} ||x_k- x||^2$. So we could get the convergence rate of the ARK algorithm:
\begin{equation}
\mathbb{E}||x_k-  x||^2 \le \lambda_{max}(C)(1-\frac{1}{\kappa(AB)^2})^{k-t0} ||x_{t0}-x||^2
\end{equation}

\subsection{Other Preconditioners}
Solving a linear system $Ax = b$ is equivalent to solving the least square problem $\sum_i ||b_i-a_i^T x||^2$. When we use a stochastic gradient descent method for this problem, a good preconditioner is the inverse Hessian matrix. And for computational efficiency, some methods have been proposed to approximate the Hessian matrix by a diagonal one, i.e. diagonal Hessian matrix. One state-of-art method is the AdaGrad method (\cite{duchi2011adaptive}).

\subsubsection{AdaGrad}
The motivation of the AdaGrad is to incorporate the geometry knowledge of the data observed in earlier iterations to adapt the weights of each dimension. 
At each step $t$, we receive a subgradient $g_t \in \partial f_t(x_t)$ of $f_t$ at $x_t$. Update $g_{1:t} = [g_{1:t-1} \quad g_t]$, $s_{t,i}= ||g_{1:t,i}||$.

Then the diagonal Hessian matrix is approximated by 
\begin{equation}
\label{adg}
H_t = \zeta I + diag(s_t).
\end{equation} 
We can use the inverse of Hessian as our preconditioner. In practice, we find that it is better to add another decay term for the inverse Hessian approximation. So, finally we use $C = \lambda_0 + H_t^{-1}$. In section 6, we compare the results of our APK algorithm and AdaSGD.

\section{The Stochastic Average Gradient based Randomized Kaczmarz Algorithm}
\label{sec:SAG}
As we introduced before, the randomized Kaczmarz algorithm can be regarded as a special form of stochastic gradient descent (\cite{needell2014stochastic}). Solving a linear system $Ax = b$ is equivalent to solving the least square problem $\sum_i f_i = \sum_i ||b_i-a_i^T x||^2$. In the randomized Kaczmarz algorithm \ref{RKalg}, the $(b_j - a_j^T,x_k) a_j$ is the gradient of a component $f_j$, and $\frac{1}{||a_j||^2}$ can be seen as the step size.
The connection between RK and SGD motivates us to bring in acceleration schemes of the SGD into the RK algorithm. 

\subsection{Stochastic Average Gradient (SAG)}
Recently, there has been a lot of work on accelerating the SGD, such as SAG (\cite{roux2012stochastic}), SDCA (\cite{shalev2013stochastic}), SVRG (\cite{johnson2013accelerating}), SAGA (\cite{defazio2014saga}). The SAG algorithm is one of the simplest of them. It requests only a small number of operations at each iteration.

As in a general stochastic optimization, the SAG method typically solve the problem of optimizing a sum of functions in this form:
$$\min_{x \in \mathbb{R}_{n}} g(x) := \frac{1}{m} \sum_{i=1}^{m} f_i(x)$$
where each $f_i$ is convex and each gradient $f'_i$ is Lipschitz continuous with constant L. To get a linear convergence rate, the average function $g(x)$ is also assumed strongly convex.
\begin{itemize}
\item $f'_i$ is Lipschitz continuous: $$||f'_i(x)-f'_i(y)|| \le L||x-y||$$
\item $g$ is strongly convex: $$g(x) \ge g(y) + g'(y)(x-y) + \frac{\mu}{2}||x-y||~2 $$
\end{itemize}

The SAG method combines the low iteration cost of the stochastic gradient descent methods with a linear convergence rate as in the full gradient methods. The method stores the most recent gradient of $f_i$(i= 1,...,m) and use the average of them to approximate the full gradient vector. At each iteration, a random training example $i_k$ is selected and $x$ is updated:
$$x^{k+1} = x^k - \frac{\alpha_k}{m}\sum_{i=1}^{m} \phi_i^k$$
where $\phi_i^k$ is updated as follows:
\begin{equation*}
\phi_i^k = \left\{
\begin{aligned} 
f'_i(x^k) \text{ if } i = i_k\\
\phi_i^{k-1} \text{ otherwise}
\end{aligned}
\right.
\end{equation*}
Let $d =\frac{1}{m}\sum_{i=1}^{m} \phi_i^k $, then at each procedure we only need to update $d$ by the following procedure:
$$d_{k+1} = d_k - \phi_i + f'_i(x^k)$$

The SAG method essentially reduced the variance between the stochastic average gradient and the full gradient (\cite{defazio2014saga}). A variance reduction approach is to use $\alpha(X-Y) + \mathbb{E}Y$ as an approximation of $\mathbb{E}X$, where $\alpha \in [0,1]$, $X$ is the SGD gradient, and $Y$ is the past stored gradient. So the variance could be changed from $Cov(X,Y)$ to $\alpha[Var(X)+Var[Y]-2Cov(X,Y)]$. SAG could be obtained from the technique by using $\alpha = 1/n$. Using the same form, when $\alpha = 1$, we could get the SAGA, which is unbiased but has larger variance than SAG. Most recently, non-uniform version of SAG, namely SAG-NUS, has also been developed to generalize SAG and SAGA, where $\alpha = \frac{1}{np}$.

SAG, SAGA, and SAG-NUS all use a constant step size in each iteration. This makes the algorithms converge fast while they are also easy to be implemented. For example, it can be demonstrated that with a constant step size of $\alpha_k = \frac{1}{2mL}$, the SAG iterations satisfy 
\begin{equation}
\mathbb{E} [||x^k - x^*||^2] \le (1-\frac{\mu}{8Lm})^k[3||x_0-x^*||^2 +C_0]
\end{equation}
where $C_0 = \frac{9 \sigma^2}{4L^2}$, and $\sigma$ is the variance of the gradient norms at the final solution $x^*$. 
Another choice $\alpha_k = \frac{1}{16L}$ has also been proved convergent in a later version of the SAG work (\cite{schmidt2013minimizing}). The SAGA method and the SAG-NUS select larger step sizes which could theoretically guarantee linear convergence. 
However, \cite{schmidt2013minimizing} indicate that in practice we could also select a larger step size for SAG. They gave two recommendations for the step size: $1/L$ and $2/(L+m\mu)$, and observed that $1/L$ always converged and performs better than the step size of $1/16L$.

\subsection{SAG-RK}
SAG is very efficient when the objective function is strongly convex. However, in a linear system, each component of the objective function $||b_i-a_i^T x||^2$ is not strongly convex if we do not add extra regularizations. In contrast, the random Kaczmarz does not have this limitation. Our second acceleration scheme combines the methods of SAG and RK. First, we use the stochastic gradient descent to make a descent direction. Then we project the point back onto the hyperplanes of each row in the linear system. 
\begin{algorithm}
\caption{The SAG-RK Algorithm}
\label{SAG-RK}
\begin{algorithmic}[1]
\FOR {$k=0,1,...$}
\STATE Select a row $j$ from $\{1,2,...m\}$ with probability $\frac{||a_j||^2}{||A||_F^2}$
\STATE Calculate the stochastic average descent $g_k$
\STATE SAG descent step: $$y_k = x_k - \alpha_k g_k$$
\STATE Project $$x_{k+1} = y_k + \frac{(b_j - a_j^T y_k)}{||a_j||^2} a_j$$
\STATE Update history information
\STATE Update $k \leftarrow k+1$
\ENDFOR
\end{algorithmic}
\end{algorithm}

The idea is motivated by the incremental constraint projection-proximal methods \cite{wang2013incremental}, which extends the projection/proximal gradient methods by using random subgradient and random constraint updates. Given a convex optimization problem $\min_{x \in X} \sum_i^N f_i(x)$ with constraints $X = \cap_{j=1}^m X_j$, the algorithm  updates $x_{k+1} = \Pi_j[x_k-\alpha_k g_i (x_k)]$ in each iteration by sampling a $j$ from the constraints and an $i$ from the components of $\sum_i^N f_i(x)$, where $g_i (x_k)$ indicates the subgradient of $f_i(x)$ at $x_k$, $\Pi_j$ denotes a Euclidean projection onto $X_j$.

In our problem we rewrite our problem as $\min_{x \in X} \sum_i ||b_i-a_i^T x||^2$ and the constraints are $X =  \cap_{j=1}^m (a_j^T x = b_j)$. For the minimum objection, we use SAG in each iteration, and then project the point onto a hyperplane $a_j^T x = b_j$. 
%

SAG-RK also has a strong relationship with SAG-NUS. From the update equation, we  get 
\begin{eqnarray*}
x_{k+1} &=& x_k  - \alpha_k g_k + \frac{(b_j - a_j^T y_k)}{||a_j||^2} a_j \\
&=& x_k  - \alpha_k g'_k 
\end{eqnarray*}
where $g'_k = g_{k-1} + \beta (- \phi_i + f'_i(x^k)) $, $\beta  = 1/m -1/\alpha_k \frac{(b_j - a_j^T y_k)}{a_j^T (y_k-y')} $, and $y'$ is the last estimate associated with $\phi_j$. So the SAG-RK is in the same family as SAG-NUS, SAG and SAGA. They only differentiate by using a different weight $\beta$.

In practice, if the step size $\alpha_k$ is small enough, we could use $x_{k+1} = y_k + \frac{(b_j - a_j^T x_k)}{||a_j||^2} a_j$ instead to update $x$ in the algorithm. This change could eliminate the calculation of gradient at $y_k$ and decrease the computational time. In fact, it is equivalent to adding a relaxation parameter to the projection step, so we call this implementation SAG-RK-relaxation. In the next section, we will compare the two implementations.

\section{Numerical Experiments}
In this section, we study the computational behavior of APK and SAG-RK, and compare them with the original RK and other acceleration schemes: the AdaGrad and the ARK.
\subsection{Synthetic Data}
We adopt two strategies of generating synthetic data which are used in RK (\cite{strohmer2009randomized}) and ARK (\cite{liu2015accelerated}) separately. For overdetermined system ($m > n$), we follow \cite{strohmer2009randomized} and let $A$ be a $m*n$ matrix whose entries are independent $N(0,1)$ random variables. In this case, the condition number of $A$ converges to:
$$\frac{\kappa(A)}{\sqrt{n}} \to \frac{1}{1-\sqrt{m/n}} $$
As $m/n$ decrease, the condition number becomes larger and larger. So, when $m=n$, we cannot control the large condition number by this method. Then we use the method in  \cite{liu2015accelerated} instead. We first generate a random $n*n$ Gaussian matrix and find its SVD: $U\Lambda V^T$. Next, we change the singular values to $\Lambda_{ii} = i^{-\alpha}$ and compute $U \Lambda V^T$ again to get a new $A$.

For the overdetermined case, we use $m = 500, n =400$; and for a square matrix, we use $m= n = 500$ and generate two matrices with $\alpha = 0.75, 0.9$ separately.

\subsection{Implementation}
\begin{itemize}
\item {\bf AdaGrad}: Using Equation (\ref{adg}) to get an approximate diagonal Hessian matrix, then we add a decay term to its inverse matrix as the preconditioner $C = \lambda_0 + H_t^{-1}$ for our linear system (In the experiments, we set $\lambda_0 =0.2$ which gets the best performance). This matrix $C$ is then employed for \ref{update} and get a AdaGrad based RK algorithm. 
\item {\bf ARK}: As we mentioned before, to implement the ARK algorithm, we have to estimate the parameter $\lambda_{min}$ first. In Liu and Wrigt's work (\cite{liu2015accelerated}), they have a strategy to approximate the real $\lambda_{min}$. Run RK for $K2$ iterations and record $x_{K2+1}$ and  $x_{K1+1}$ where $K1 = max(1,K2-10m)$. Based on the convergence rate, they estimate the $\lambda_{min}$ as follows\footnote{In the original paper, they normalize the matrix $A$ so that $||A||_F^2 = m$}.:
$$||A||_F^2 \left[ 1-\left( \frac{||A x_{K2}-b||}{||A x_{K1}-b||}\right)^{\frac{0.5}{K_2-K_1}}\right]$$
\end{itemize}
So the ARK needs a long burn-in time to determine a good parameter. In this paper we fix $K2 = 15m$.

\begin{figure}[ht]
	\centering
	\subfigure[]{
		\includegraphics[width=0.45\linewidth]{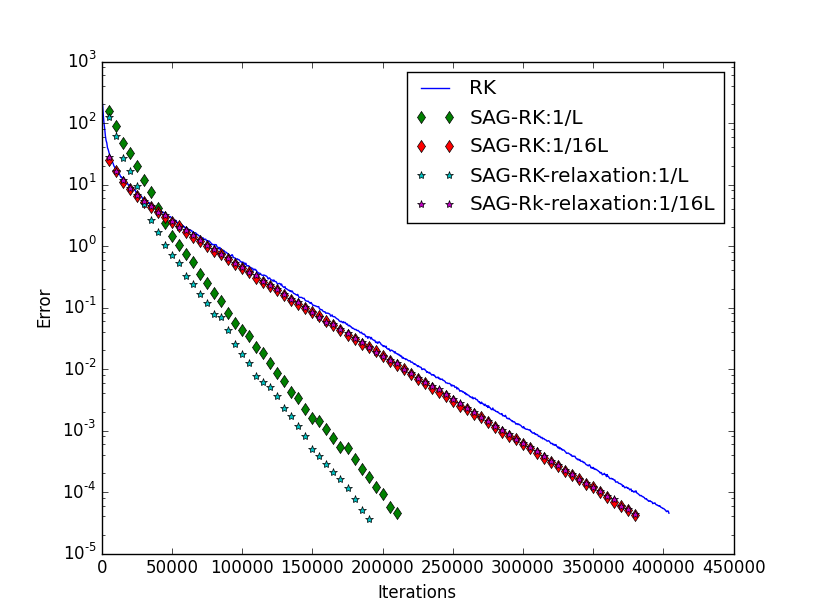}
		\label{fig:stepsize}
	}
	\subfigure[]
	{
		\includegraphics[width=0.45\linewidth]{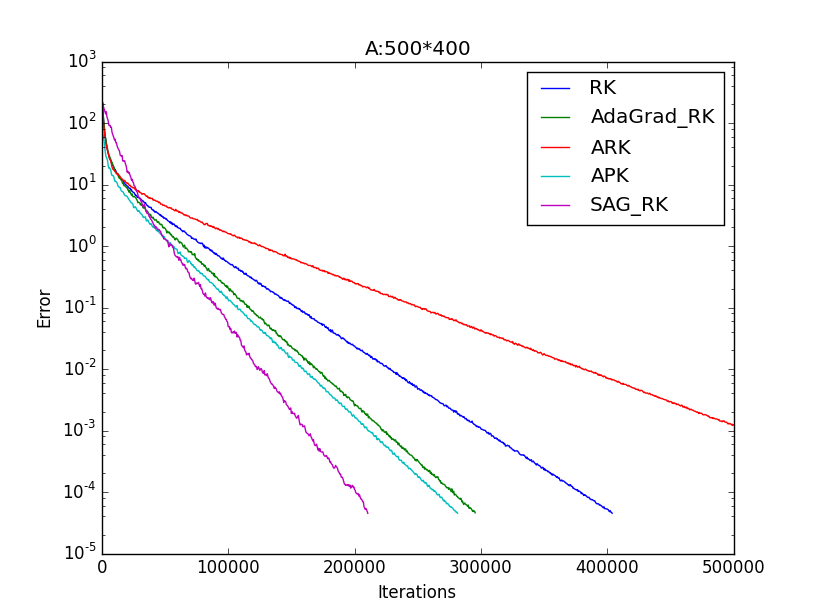}
		\label{fig1}
	}
	\subfigure[]
	{
		\includegraphics[width=0.45\linewidth]{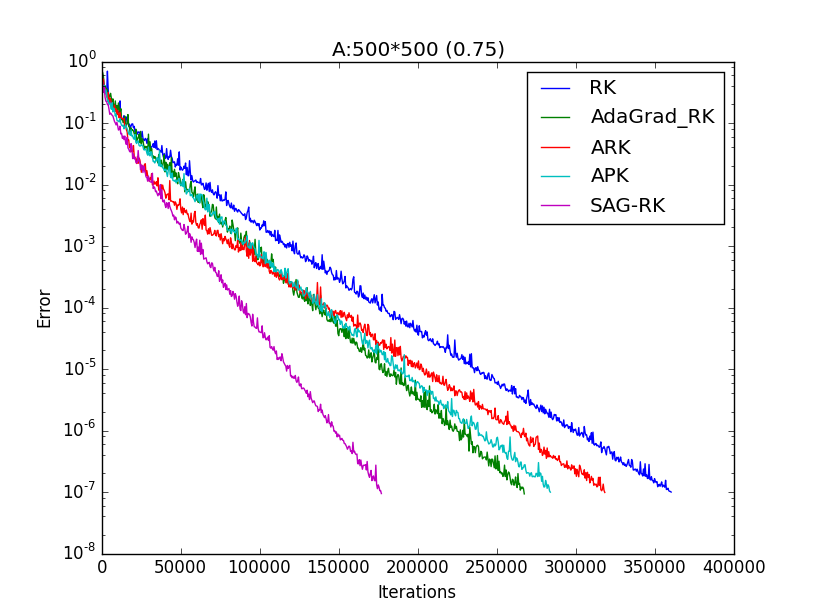}
		\label{fig2}
	}
	\subfigure[]
	{
		\includegraphics[width=0.45\linewidth]{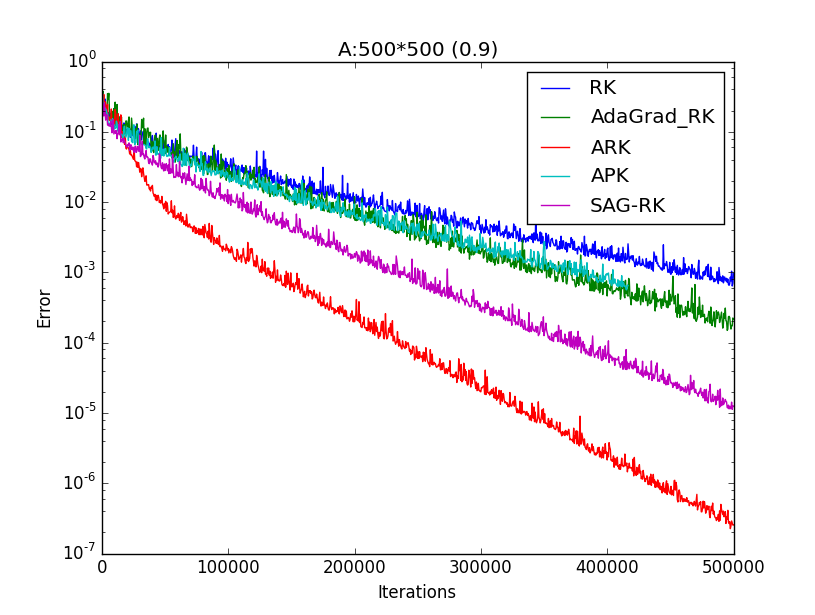}
		\label{fig3}
	}
	\caption{Numeric Experiments. (a): Comparison of step sizes in SAG-RK. (b): $A_1 \in \mathbb{R}^{500*400}$ and $\kappa (A_1) = 180.7$ (c): $A_2 \in \mathbb{R}^{500*500}$, and $\kappa (A_2) = 167.9$ (d): $A_3 \in \mathbb{R}^{500*500}$, and $A$ has a larger condition number  $\kappa (A_3) = 367.6$}
\end{figure}

For our own two algorithms, we follow the procedures that are introduced in Section\ref{sec:APK} and Section \ref{sec:SAG}. Selecting a proper step size is important for SAG, so we try different choices and compare them first. In Figure \ref{fig:stepsize}, we show residual errors for SAG-RK and its relaxation implementation with two stepsizes separately on a $500*400$ matrix. At each stepsize, the two implementations have almost the same performance, and a larger step size leads to significantly faster convergence\footnote{It has a different story when $m>>n$, but it is beyond the scope of this paper, because we only consider the ill-conditioned case.}.

Then we choose the step size of SAG-RK as $1/L$, and show the performance of all algorithms in Figure \ref{fig1}, \ref{fig2}, \ref{fig3}.  First, we compare the APK and AdaGrad, the APK has only slightly better convergence rates over AdaGrad. However, AdaGrad needs to update the diagonal matrix at each iteration and costs too much for the RK algorithm, while APK only update the preconditioner after a long interval. So APK is much more computationally efficient.

Then we compare our algorithm with the ARK algorithm. The ARK algorithm does not perform well when the condition number is not large enough, it is even worse than the RK algorithm. But it has the best performance when we have the largest condition number. This agrees with the conclusions in \cite{liu2015accelerated} that the ARK only suits to seriously ill-conditioned problems. In contrast to ARK, our algorithms performs consistently better than RK. In particular, the SAG-RK performs best on $A_1$ and $A_2$ and second best on $A_3$.

At last, we show how our algorithms improve RK on computational time in Table \ref{tb:time}. As SAG-RK can be seen a combination of SAG and RK, we compare it with the two algorithms. As we discussed before, to save the operation time at each iteration, we use another implementation, namely SAG-RK-relaxation (In the table, we call it SAG-RK2 to save space). All algorithms check the residual error every $10m$ iterations and stop at $||b-Ax||/||b||<10^{-7}$.

From Table \ref{tb:time}, we could see that SAG-RK has a much shorter computational time than RK and SAG. Furthermore, the SAG-RK-relaxation indeed improves the computational efficiency and outperforms all others. APK also has a descent performance, but it cannot get significant improvement over RK in all cases.
\begin{table}[tb]
\caption{Computational time (seconds) until $||b-Ax||/||b||<10^{-7}$} \label{tb:time}
\begin{center}
\begin{tabular}{llll}
{\bf Model}  &{\bf $A_1$} &{\bf $A_2$}	 &{\bf $A_3$}\\ 
\hline \\
RK         &50.80	&43.42	&159.17 \\
SAG             &339.29 	&48.25	&182.65\\
SAG-RK             &36.87	&31.35	&116.13\\
SAG-RK2             &31.95 		&26.87	 &99.63\\
APK             &38.41		&40.43	 &156.86\\
\end{tabular}
\end{center}
\end{table}

\section{Conclusion}
In this paper, we propose two methods to accelerate the randomized Kaczmarz algorithm, namely APK and SAG-RK. They both take advantage of the history information in past iterations. APK use past estimates of $x$ to get an approximate right preconditioner for the linear system, while SAG-RK use past gradients to get an approximate full gradient and combine the SAG step with a Kaczmarz projection. The APK provides a new vision to develop preconditioners based on history information. And the performance of SAG-RK exceeds both SAG and RK consistently. Future work includes extension to the inconsistent linear systems as well as the sparse case, and theoretic demonstration of the convergence rate.

\bibliographystyle{plain}
\bibliography{Kaczmarz}

\end{document}